\documentclass[11pt]{article}
\usepackage[margin=1in]{geometry}
\usepackage{amsmath,amssymb,amsthm}
\usepackage{graphicx}
\usepackage{booktabs}
\usepackage[hidelinks]{hyperref}

\newtheorem{theorem}{Theorem}

\title{Small $5$-chromatic unit-distance graphs in disks of radius $3/2$, $5/4$ and $23/20$}
\author{Chandan Shukla\\ \small\texttt{shukla.chandan12@gmail.com}}
\date{November 2026}

\begin{document}
\maketitle

\begin{abstract}
We give three explicit unit-distance graphs in the plane with chromatic number~$5$
whose vertices lie in small disks: a graph with $773$ vertices inside an open disk of
radius $3/2$, one with $5713$ vertices inside an open disk of radius $5/4$, and one with
$8168$ vertices inside an open disk of radius $23/20$. The smallest examples previously
reported for these radii have $11725$, $11725$ and $31171$ vertices.
All coordinates are exact elements of $\mathbb{Q}(\sqrt3,\sqrt5,\sqrt{11})$.
In each case $5$-chromaticity is certified by an explicit $5$-colouring and a
DRAT proof that no $4$-colouring exists. The $773$-vertex graph is vertex-critical.
\end{abstract}

\section{Introduction}

A unit-distance graph is a graph whose vertices are points of the plane, two of them
adjacent when they are at distance~$1$. The chromatic number of the plane is the
largest chromatic number of such a graph. De~Grey~\cite{deGrey} found a
unit-distance graph with chromatic number~$5$, and a second proof was given by
Exoo and Ismailescu~\cite{EI}. The search for small examples went on in the
Polymath16 project; the smallest known $5$-chromatic unit-distance graph has $509$
vertices and is due to Parts~\cite{Parts}, after a $553$-vertex graph of
Heule~\cite{Heule}. Parts also gave a proof of $\chi(\mathbb{R}^2)\ge5$ that can be
checked by hand~\cite{PartsHuman}, and de~Grey and Parts~\cite{dGPorder} gave lower
bounds on the number of vertices of $k$-chromatic unit-distance graphs. Related
colouring questions for the plane are studied in~\cite{PartsPercent,dGPtiling,
PartsMultifold,PartsInterval,PartsPlane}, for two forbidden distances
in~\cite{PartsTwo}, for odd distances in~\cite{PartsOdd}, and in higher dimension
in~\cite{PartsR8}.

These graphs are spread out. Parts' $509$-vertex graph has diameter $4.950$, and the
smallest disk containing it has radius $2.524$. The natural local question,
studied in~\cite{PolyDisk}, is how small a disk can carry a $5$-chromatic
unit-distance graph. Let $r_5$ be the infimum of the radii $r$ such that some finite
unit-distance graph with all vertices in a closed disk of radius~$r$ is not
$4$-colourable. The Polymath paper~\cite{PolyDisk} shows
\[
  0.86947 < r_5 \le 1 .
\]
The upper bound comes from a graph with $31171$ vertices and $180942$ edges, made of a
pair of graphs of type $4H4$ placed on a disk of radius exactly~$1$; the authors note
that they did not try to reduce its size. In the Polymath16 discussion~\cite{Thread17}
Parts also reported $5$-chromatic graphs on disks of radius $1.25$ ($11725$ vertices),
$1.09$ ($75793$ vertices) and $1.04$ ($311509$ vertices).

Once $r\ge1$ the question is no longer whether a $5$-chromatic graph fits in a disk of
radius~$r$, but how few vertices it needs. Since a graph in a disk of radius~$1$ also
lies in every larger disk, the previously known minimum for radius $3/2$ and $5/4$ is
$11725$, and for radius $23/20=1.15$ it is $31171$.

\begin{theorem}\label{thm:main}
There are unit-distance graphs $G_{773}$, $G_{5713}$ and $G_{8168}$ with chromatic
number~$5$ and with the numbers of vertices and edges and the enclosing disks shown in
Table~\ref{tab:main}. All vertices lie in the open disk named in the table.
The graph $G_{773}$ is vertex-critical: deleting any one vertex leaves a
$4$-colourable graph.
\end{theorem}

\begin{table}[ht]
\centering
\begin{tabular}{lrrlccr}
\toprule
graph & vertices & edges & disk (centre; radius) & smallest radius & diameter & previous\\
\midrule
$G_{773}$  &  773 &  3817 & $(0,0);\ 3/2$                 & 1.49395 & 2.98790 & 11725\\
$G_{5713}$ & 5713 & 25482 & $(0,0);\ 5/4$                 & 1.24821 & 2.49641 & 11725\\
$G_{8168}$ & 8168 & 41556 & $(7/16,-1693/5000);\ 23/20$    & 1.14962 & 2.29747 & 31171\\
\bottomrule
\end{tabular}
\caption{The three graphs. ``Smallest radius'' is the radius of the smallest disk
containing the vertex set, rounded up; ``previous'' is the smallest number of
vertices previously reported for a $5$-chromatic unit-distance graph in a disk of the
same radius.}
\label{tab:main}
\end{table}

The theorem does not change the bounds on $r_5$. What it shows is that the radius
constraint is much cheaper than the earlier examples suggested: at radius $3/2$ we
need about $1.5$ times as many vertices as the unrestricted record of $509$, where the
earlier examples needed more than twenty times as many.

\section{The graphs}

All vertices have coordinates in the field $K=\mathbb{Q}(\sqrt3,\sqrt5,\sqrt{11})$,
which has degree $8$ over $\mathbb{Q}$ with basis
\[
  (\beta_0,\dots,\beta_7)=(1,\sqrt3,\sqrt5,\sqrt{15},\sqrt{11},\sqrt{33},\sqrt{55},\sqrt{165}).
\]
A vertex is stored as $16$ integers $a_0,\dots,a_{15}$ and stands for the point
\[
  p=\frac{1}{96}\Bigl(\sum_{k=0}^{7}a_k\beta_k,\ \sum_{k=0}^{7}a_{8+k}\beta_k\Bigr).
\]
For two vertices $p,q$ the number $96^2|p-q|^2$ is an element of $K$ with integer
coordinates in this basis, computed exactly from the products
$\beta_i\beta_j$. The pair is an edge if and only if this element equals $96^2$, that is,
its rational part is $9216$ and its other seven coordinates are $0$. Every graph here
is the full unit-distance graph on its vertex set, so the coordinates determine the
edges.

The coordinates of all three graphs, their $5$-colourings and a checking program are
in the ancillary files. Some further data:

\begin{center}
\begin{tabular}{lrrrr}
\toprule
graph & min degree & max degree & triangles & edge directions (up to sign)\\
\midrule
$G_{773}$  & 4 &  42 &  1335 & 36\\
$G_{5713}$ & 4 &  66 &  9102 & 58\\
$G_{8168}$ & 4 & 104 & 12022 & 58\\
\bottomrule
\end{tabular}
\end{center}

The vertex of $G_{773}$ at the origin has degree $42$. The five colour classes in the
colouring we give have sizes $181$, $141$, $133$, $160$ and $158$.

\begin{figure}[ht]
\centering
\includegraphics[width=0.62\textwidth]{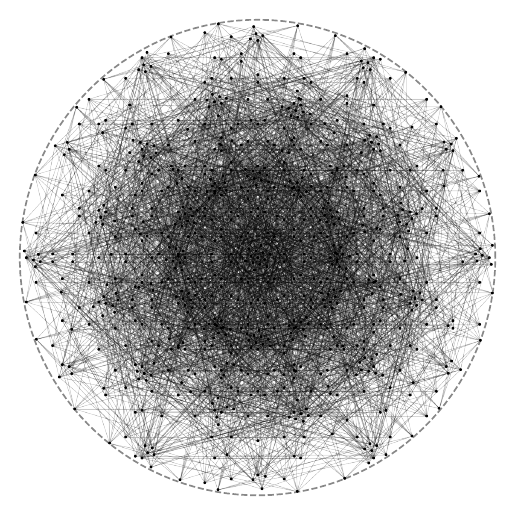}
\caption{$G_{773}$: all $773$ vertices and $3817$ edges. The dashed circle has
radius $3/2$ and centre $(0,0)$.}
\label{fig:773}
\end{figure}

\begin{figure}[ht]
\centering
\includegraphics[width=\textwidth]{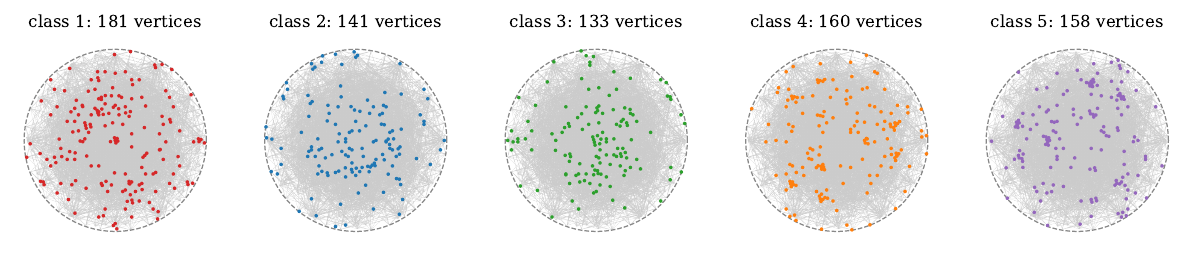}
\caption{The five colour classes of a proper $5$-colouring of $G_{773}$
(edges in grey).}
\label{fig:classes}
\end{figure}

\begin{figure}[ht]
\centering
\includegraphics[width=0.46\textwidth]{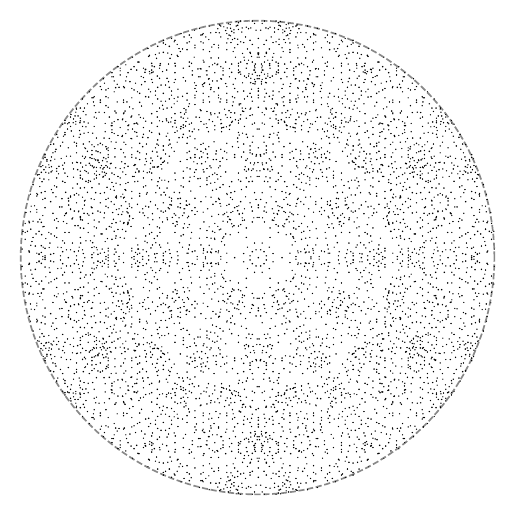}\hfill
\includegraphics[width=0.46\textwidth]{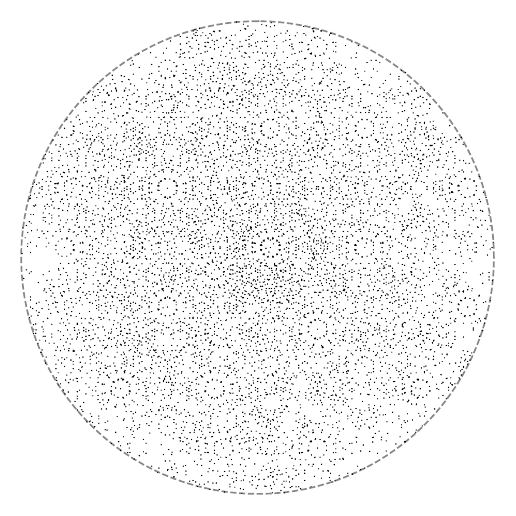}
\caption{The vertex sets of $G_{5713}$ (left, circle of radius $5/4$) and
$G_{8168}$ (right, circle of radius $23/20$ centred at $(7/16,-1693/5000)$).
The edges are omitted; there are too many to draw legibly.}
\label{fig:big}
\end{figure}

\section{Verification}

For each graph we checked the following.

\begin{enumerate}
\item \emph{Geometry.} The points are distinct, and every point lies strictly inside
the disk of Table~\ref{tab:main}. For the second check each $\sqrt{m}$ is replaced by
rational lower and upper bounds with error below $10^{-30}$, which gives a rational
upper bound for $|p-c|^2$; this bound is below $r^2$ for every vertex.
\item \emph{Edges.} Candidate pairs are those whose floating-point distance is within
$10^{-6}$ of~$1$; each candidate is then decided in exact integer arithmetic as
described above. Pairs outside this window are certainly not at distance~$1$.
\item \emph{Upper bound.} The listed colouring uses $5$ colours and no edge joins two
vertices of the same colour.
\item \emph{Lower bound.} The graph is not $4$-colourable. We encode $4$-colourability
in the usual way, with one Boolean variable for each vertex and colour, and fix the
colours of the three vertices of one triangle to $1,2,3$; this loses no generality,
since colours can be permuted. The formula was shown unsatisfiable by
Kissat~\cite{Kissat}, and the resulting DRAT proof was checked with
\texttt{drat-trim}~\cite{DratTrim}.
\end{enumerate}

For $G_{773}$ the formula has $3092$ variables and $16044$ clauses; the proof was
found in about $17$ seconds and checked in about $15$ seconds. For $G_{5713}$
($22852$ variables, $107644$ clauses) Kissat needed about $35$ minutes and the proof
($5.5$\,GB in text form) was checked in about $31$ minutes. For $G_{8168}$
($32672$ variables, $174395$ clauses) Kissat needed about $6$ minutes and the check
took about $6$ minutes. All runs used a single core of a laptop. We also checked that $G_{773}$ is vertex-critical: for each of its $773$ vertices
$v$ we found an explicit $4$-colouring of $G_{773}-v$ and checked that it is proper.

The graphs $G_{5713}$ and $G_{8168}$ have minimum degree $4$, but we have not tried to
make them critical, and they can very likely be reduced.

\section{Remarks}

The main open question about disks remains the one from~\cite{PolyDisk}: is there a
$5$-chromatic unit-distance graph in a disk of radius less than~$1$? A disk of radius
less than $1$ contains no unit hexagon, and all small known $5$-chromatic graphs use
many of them, so a different kind of construction seems to be needed.

For radii between $1$ and $3/2$ the vertex counts are far from settled. The graphs
here suggest that the minimum number of vertices drops quickly as the radius grows
from about $1.15$ to $1.5$; it would be interesting to know how it behaves near
radius~$1$, where the only known example still has $31171$ vertices.

\section*{Ancillary files}

\begin{itemize}
\item \texttt{g773.txt}, \texttt{g5713.txt}, \texttt{g8168.txt}: one vertex per line,
$16$ integers $a_0,\dots,a_{15}$ as above.
\item \texttt{g773\_5colouring.txt}, \texttt{g5713\_5colouring.txt},
\texttt{g8168\_5colouring.txt}: one colour ($0$--$4$) per line.
\item \texttt{verify.py}: checks the geometry, recomputes all edges exactly, checks a
colouring, and writes the $4$- and $5$-colouring formulas. Example:
\texttt{python3 verify.py g773.txt 0 0 3/2 --colouring g773\_5colouring.txt --cnf g773\_}
\end{itemize}

\section*{Acknowledgements}

Some of the computer searches and checks were carried out with the help of AI
coding assistants. All results stated here were verified as described above.

\end{document}